\documentclass[11pt]{article}

\usepackage{latexsym,amssymb}
\usepackage{amsmath,amsthm}
\usepackage{fancyhdr}
\usepackage{geometry}
\usepackage{fourier}

\usepackage[colorlinks=true, linkcolor=DarkSlateBlue, %
	citecolor=purple, urlcolor=blue, pdfauthor={James Montaldi}]{hyperref}

\usepackage[backend=biber]{biblatex}
\usepackage[svgnames]{xcolor}
\usepackage{pstricks,pst-plot}

\usepackage{caption}
\DeclareCaptionLabelFormat{nocolon}{\textsc{#1} #2}
\renewcommand{\sectionmark}[1]{%
\markboth{\thesection.\ #1}{}}

\numberwithin{equation}{section}
\numberwithin{figure}{section}
\numberwithin{table}{section}
\newtheorem{theorem}{Theorem}[section]
 \newtheorem{lemma}[theorem]{Lemma}
 \newtheorem{proposition}[theorem]{Proposition}
 \newtheorem{corollary}[theorem]{Corollary}

\theoremstyle{definition}
 \newtheorem{definition}[theorem]{Definition}
 \newtheorem{example}[theorem]{Example}
 
 \newtheorem{remark}[theorem]{Remark}
 \newtheorem{GR}[theorem]{Geometric remark}

\newcommand{\defn}[1]{{\bfseries\itshape #1}}

\newcommand{\dbyd}[2]{\frac{\partial #1}{\partial #2}}

\def\restr#1{\,\rule[-1ex]{0.4pt}{2ex}\lower1.0ex\hbox{\scriptsize $\;#1$}}

\def\Reeb{\mathcal{R}}
\def\half{\tfrac12}
\def\ij{i\kern-1pt j}
\def\Image{\mathrm{Image}}

\def\d{\mathrm{d}}
\def\D{\mathrm{D}}
\def\df{\mathrm{d}\!f}
\def\dh{\mathrm{d}\!h}
\def\dH{\mathrm{d}\!H}
\def\dpp{\mathrm{d}\!p} % \dp is a TeX primitive
\def\dq{\mathrm{d}\!q}
\def\dt{\mathrm{d}\!t}
\def\dx{\mathrm{d}\!x}

\def\dz{\mathrm{d}\!z}

\def\aa{\mathbf{a}}
\def\bb{\mathbf{b}}

\def\uu{\mathbf{u}}
\def\vv{\mathbf{v}}
\def\ww{\mathbf{w}}
\def\xx{\mathbf{x}}
\def\yy{\mathbf{y}}
\def\RR{\mathbb{R}}

\def\Hess{{\mathrm{Hess}}}
\def\tr{\mathop\mathrm{tr}}

\def\zero{\mathbf{0}}

\def\kK{\mathcal{K}}
\def\kX{\mathcal{X}}

\begin{document}
\thispagestyle{empty}

\noindent{\huge\textbf{Equilibria and bifurcations in \\[4pt]
contact dynamics}}

\bigskip
\bigskip

{\color{Gray}
\noindent{\Large\textbf{James Montaldi}}

\medskip

\noindent\textbf{University of Manchester}

\noindent\today

\medskip
}

{\small
%%%%%%%%%%%%%%%%%%%%%%%%%%%
\noindent\hrulefill % start of abstract

\smallskip

\noindent{\large\scshape\color{Gray} Abstract}

\medskip

\noindent 
We provide a systematic study of equilibria of contact vector fields and the bifurcations that occur generically in 1-parameter families, and express the conclusions in terms of the Hamiltonian functions that generate the vector fields. 

Equilibria occur at points where the zero-level set of the Hamiltonian function is either singular or is tangent to the contact structure.  The eigenvalues at an equilibrium have an interesting structure: there is always one particular real eigenvalue of any equilibrium, related to the contact structure, that we call the principal coefficient, while the other eigenvalues arise in quadruplets, similar to the symplectic case except they are translated by a real number equal to half the principal coefficient.  

There are two types of codimension 1 equilibria, named Type I, arising where the zero-set of the Hamiltonian is singular, and Type II where it is not, but there is a degeneracy related again to the principal coefficient and the contact of the zero level-set of the Hamiltonian with the contact structure.  Both give rise generically to saddle-node bifurcations.  

Some special features include: (i) for Type II singularities, Hopf bifurcations cannot occur in dimension 3, but they may in dimension 5 or more; (ii) for Type I singularities, a fold-Hopf bifurcation can occur with codimension 1 in any dimension, and (iii) again for Type I, and in dimension at least 5, a fold-multi-Hopf bifurcation (where several pairs of eigenvalues pass through the imaginary axis simultaneously together with one through the origin) may also occur with codimension 1. 

\medskip

\noindent \emph{MSC 2020}: 37G10; 37J55; 53E50;   \\[6pt]
\noindent \emph{Keywords}:  bifurcations, contact structure,  Hamiltonian system,  

\noindent\hrulefill % end of abstract
%%%%%%%%%%%%%%%%%%%%%%%%%%%
}

\bigskip

\tableofcontents

%%%%%%%%%%%%%%%%%%%%%%%%%%%
%%%%%%%%%%%%%%%%%%%%%%%%%%%
\section*{Introduction} \sectionmark{Introduction}

A contact vector field on a contact manifold is one whose flow preserves the contact structure. We study the most elementary aspects of the dynamics of such vector fields; namely, equilibria, their stability and their generic 1-parameter bifurcations.  

There has been considerable interest in contact vector fields in recent years, in several different directions.  For example, they play a role in  thermodynamics (see for example A.~Bravetti \cite{Bravetti-ContactGeom+Thermo} and references therein as well as D.~Gromov \cite{Gromov2016}), in Hamiltonian-like systems with dissipation, both classical \cite{deLeon+LainzValcazar} and quantum \cite{Marmo-etal-2018}, in fluid dynamics \cite{Ghrist-Fluids, Grmela+Ottinger}, and others.  An interesting application of contact geometry to neuroscience can be found in an article of  Petitot \cite{Petitot}.
For more details, examples and discussions see the review by Bravetti \cite{Bravetti2017}.

There are few dynamical studies beyond setting up a model, though Gromov and Cairnes \cite{Gromov-Caines2015} consider dynamics for a diatomic gas, and Liu et al \cite{Liu-etal-2018} consider periodic motion in some restricted settings, and a recent paper of Entov and Polterovich \cite{Entov+Polterovich-2023} discuss trajectories with special properties. It was also found by Bravetti and Tapias \cite{Bravetti+Tapias2015}  that there is an invariant measure defined on the open dense subset of the open submanifold where the Hamiltonian is non-zero, and Bravetti et al.\,\cite{Bravetti-deLeon+}  consider the type of dynamics on the complement, that is where $H=0$.

A contact structure $\xi$ on a manifold $M$ consists of, for each $x\in M$, a hyperplane $\xi(x)\subset T_xM$ such that this hyperplane field is maximally non-integrable. The easiest way to define the maximally non-integrable property is to choose any (local) 1-form $\eta$ such that $\ker\eta(x)=\xi(x)$ in the domain of $\eta$ (these are called contact 1-forms), and the non-integrability requirement is that the volume form $\eta\wedge(\d\eta)^n$ should not vanish anywhere (this is independent of the choice of contact 1-form $\eta$).  For general background on contact structures, the reader can consult \cite{A+M-FoM,Arnold,Arnold-Givental,Geiges,Libermann-Marle}. 

In many areas, for example thermodynamics and jet bundles, the form $\eta$ plays a primary role, with $\xi$ being a secondary construction.  On the other hand, many authors put $\xi$ in the forefront, and choose a contact form for calculational convenience.    

A further approach is taken by Grabowska and Grabowski \cite{Gra-Gra}; they explicitly put $\xi$ at the forefront by considering the line bundle $\xi^\circ\subset T^*M$ whose fibre is the annihilator of $\xi$.  Any contact form $\eta$ for $\xi$ is a section of this bundle. Hamiltonians are then defined to be functions on this line bundle that are homogeneous of degree 1.  In other words, they consider all possible contact 1-forms together.

We begin the paper by recalling the basic well-known properties of contact geometry and contact vector fields, and describe how every contact vector field has a unique Hamiltonian function which generates it.  In \S\ref{sec:equilibria} we begin the study of equilibria: what are the conditions on the Hamiltonian for a point to be an equilibrium, and what are the conditions for an equilibrium to be nondegenerate?  It turns out that generically equilibria occur where the zero-level of the Hamiltonian is tangent to the contact hyperplane.  A central result (Theorem \ref{thm:linear} and its corollary) states that at an equilibrium one of the eigenvalues, which we call the principal coefficient,  is real and related to the Reeb vector field, and the others arise in quadruplets similar to the symplectic case, but translated by a real number equal to one half the principal coefficient.  We also show how Hopf bifurcations can arise in dimension at least 5.  We identify two ways in which an equilibrium can degenerate. The first, which we call Type I degeneracy, arises where $H$ has a critical point on its zero-level set, and the second, Type II, when the restriction to the contact plane is degenerate.  

The following two sections describe the nature of the Type I and Type II degeneracies, respectively and a brief analysis of the resulting saddle node bifurcations. 

We end with a short discussion of Legendre vector fields; that is, vector fields on Legendre submanifolds. In particular we show that, given any vector field on a Legendre submanifold, there is an extension of it to a contact vector field on the ambient contact manifold. A consequence of this is that the bifurcation theory of Legendre vector fields is the same as for ordinary vector fields in that dimension. 

The paper ends with a short appendix containing an elementary singularity theoretic calculation for recognizing folds and their versal unfoldings.

%%%%%%%%%%%%%%%%%%%%%%%%%%%
%%%%%%%%%%%%%%%%%%%%%%%%%%%
\section{Background}
\label{sec:background}

Here we establish some notation used throughout.  Every object we consider will be assumed to be smooth.  Let $(M,\xi)$ be a contact manifold, with $\dim M = 2n+1$. This means that $\xi$ is a subbundle of $TM$ of rank $2n$ which is maximally non-integrable. 

The non-integrability condition is most easily expressed in terms of 1-forms vanishing on $\xi$: let $\eta$ be any 1-form on $M$  satisfying $\ker\eta=\xi$ (possibly defined locally).  Then the non-integrability condition states that $\eta\wedge(\d\eta)^n$ is nowhere zero. For details, see for example \cite{A+M-FoM, Arnold,Arnold-Givental,Geiges, Libermann-Marle}. 
We say that a 1-form $\eta$ with the property that at each point of its domain of definition, $\ker\eta=\xi$ is a \defn{contact form} for $\xi$. Such an $\eta$ is determined up to non-zero scalar multiples: that is, if $\eta_1,\eta_2$ are contact forms for $\xi$ then there is a nowhere zero function $f$ on $M$ such that $\eta_2=f\eta_1$. 

We write $\kX_\xi$ for the space of vector fields on $M$ whose flow preserves the contact structure.
Let $X\in\kX_\xi$.  

The easiest way of determining whether a given vector field $X$ on $M$ preserves $\xi$ is to introduce a contact 1-form $\eta$.  One sees that $X$ preserves $\xi$ if and only there is a function $f$ (possibly zero) such that $L_X\eta=f\eta$.  

\begin{definition}\label{def:Reeb VF}
Given a contact 1-form $\eta$ on $(M,\xi)$ the \defn{Reeb vector field} $\Reeb$ is uniquely determined by the conditions  
$$\iota_\Reeb\eta = 1\quad\text{and} \quad\iota_\Reeb\d\eta=0.$$
\end{definition}

Note that the Reeb vector field is dependent on the choice of contact form; indeed not even its direction is intrinsically associated to $\xi$. Introducing a contact form $\eta$ allows one to state the following well-known criterion.

\medskip

\begin{proposition}\label{prop:contact VF}
Let $Y$ be any (smooth) vector field on $M$, let $\eta$ be any contact 1-form for $(M,\xi)$ and let $h=-\eta(Y)$. Then $Y\in\kX_\xi$ if and only if 
$$\iota_Y\d\eta = \d h-\Reeb(h)\eta,$$
where $\Reeb$ is the Reeb vector field associated to $\eta$ and $\Reeb(h)$ is the derivation of\/ $h$ along the vector field $\Reeb$; that is, $\Reeb(h)(x)=\dh_x(\Reeb(x))$. 
\end{proposition}

\medskip

\begin{proof}
Let $Y$ be a vector field on $M$.  We have 
$L_Y\eta=\d(\iota_Y\eta)+ \iota_Y\d\eta = -\dh+\iota_Y\d\eta$.
That is, for any vector field $Y$
$$\iota_Y\d\eta = \dh +L_Y\eta.$$

Suppose first that $Y$ preserves $\xi$; that is, $L_Y\eta=f\eta$ for some function $f$. Now $\iota_\Reeb\d\eta=0$,  and hence
\begin{eqnarray*}
0&=&\iota_\Reeb\iota_Y\d\eta \\
&=& \iota_\Reeb(\dh +f\eta)\\
&=& \Reeb(h) + f
\end{eqnarray*}
whence $f=-\Reeb(h)$. That is, any vector field $Y\in\kX_\xi$ satisfies $L_Y\eta=-\Reeb(h)\eta$, and hence the expression for $\iota_Y\d\eta$ follows. 

Conversely, if $\iota_Y\d\eta=\dh-\Reeb(h)\eta$ then $L_Y\eta=-\dh+\dh-\Reeb(h)\eta=-\Reeb(h)\eta$ so that $Y\in\kX_\xi$. 
\end{proof}

\begin{definition}\label{def:Hamiltonian}
Let $X$ be a contact vector field on $(M,\xi)$ and let $\eta$ be a contact 1-form. The function
$H=-\iota_X\eta$ is called the \defn{Hamiltonian} of the vector field (associated to $\eta$).  
\end{definition}
We will usually take $\eta$ as given, but if $\eta$ were replaced by $f\eta$ for some non-zero smooth function $f$, then $H$ would be replaced by $fH$.

By Proposition\,\ref{prop:contact VF}, the Hamiltonian satisfies
\begin{equation}\label{eq:Ham VF}
\iota_X\d\eta = \dH-\Reeb(H)\eta.
\end{equation}

The definition gives a linear map $\kX_\xi\to C^\infty(M,\RR)$, $X\mapsto-\iota_X\eta$. This is in fact an isomorphism, whose inverse is as follows. 

Given a `Hamiltonian' function $H$, the associated vector field $X=X_H$ is defined implicitly by the equations 
\begin{subequations} \label{eq:Hamiltonian->VF}
\begin{eqnarray}
\eta(X_H) &=& -H \label{eq:normal}\\
\iota_{X_H}\d\eta &=& \dH-\Reeb(H)\,\eta   \label{eq:tgt}
\end{eqnarray}
\end{subequations}
The first equation determines the normal component of $X_H$, and the second the `tangential' component (ie the component on $\xi$). Proposition\,\ref{prop:contact VF} shows that such vector fields preserve the contact structure.  

The fact that this is an isomorphism allows us to parametrize the space of smooth contact vector fields by using smooth functions on $M$. Thus one likes to describe any property of the vector field in terms of its Hamiltonian function.  For example, from the definition of $H$, it follows that $X_H(x_0)\in\xi(x_0)$ if and only if $H(x_0)=0$.  

Note that $X_{H+C} = X_H-C\Reeb$, so adding a constant to $H$ changes the dynamics. In particular the Reeb vector field itself is the contact vector field associated to the constant function $H=-1$. 

One important contrast with Hamiltonian vector fields on a symplectic manifold is that in general the Hamiltonian is not a conserved quantity.  Indeed, applying \eqref{eq:tgt} to $X_H$ gives $0=\dH(X_H)-\Reeb(H)\eta(X_H)$ which by \eqref{eq:normal} leads to
\begin{equation}\label{eq:derivative of H}
\frac{\d}{\dt}H = \dH(X_H) =  -\Reeb(H)H.
\end{equation}
In particular only $H^{-1}(0)$ is an invariant level set in general.  This formula for $\frac{\d}{\dt} H(t)$ shows also that $H^{-1}(0)$ is attracting if and only if $\Reeb(H)>0$ along this hypersurface; this is important when using contact vector fields to model dissipation.   

Another useful property of the Reeb vector field is that it determines whether the form $\eta$ is preserved by $X$: using Cartan's formula, one shows that if $X\in\kX_\xi$ then 
$$L_{X}\eta = -\Reeb(H)\eta$$
where $H=-\eta(X)$. In particular $X$ preserves $\eta$ if and only if $\Reeb(H)=0$, which in turn is equivalent, by \eqref{eq:derivative of H}, to $X_H$ preserving every level set of the Hamiltonian; such vector fields are often called \emph{strict} contact, or \emph{conservative}, vector fields. Moreover, if $\nu=\eta\wedge(\d\eta)^n$ is the contact volume form, then it follows from the above that
$$L_X\nu = -(n+1)\Reeb(H)\nu.$$

\paragraph{Darboux coordinates}
Recall (for a proof see for example \cite{Geiges,Libermann-Marle}) that if $x_0\in M$ then there is a neighbourhood of $x_0$ and coordinates $q_1,\dots,q_n,p_1,\dots,p_n,z$ such that 
\begin{equation}\label{eq:Darboux}
\eta = \dz - p_i\dq_i
\end{equation}
(where the summation convention is understood).  These are called \defn{Darboux} or \defn{canonical coordinates}. For this $\eta$, the contact hyperplane $\xi$ has basis
\begin{equation}\label{eq:basis for xi}
\left\{\partial_{p_j},\; \partial_{q_j}+p_j\partial_z\right\}\quad (j=1,\dots,n).
\end{equation}
Notice that on $\xi$ this basis is \emph{canonical}, in the sense that 
$$\d\eta\,(\partial_{q_j}+p_j\partial_z, \: \partial_{p_i}) = \delta_{\ij},\quad\text{etc.},$$
where $\delta_{\ij}$ is the Kronecker delta. 

On $\RR^{2n+1}$ with canonical/Darboux coordinates as above, we have $\Reeb=\partial_z$. 

To describe the relation between $H$ and $X_H$ one can use the method of coefficients. Given a (Hamiltonian)  function $H$, write $X_H=a_i\partial_{q_i}+b_i\partial_{p_i}+c\partial_z$. Then equations \eqref{eq:Hamiltonian->VF} give
$$c-p_ia_i=-H,\quad \text{and}\quad a_i\dpp_i -b_i\dq_i = (H_{q_i}\dq_i+H_{p_i}\dpp_i+H_z\dz)-H_z(\dz-p_i\dq_i).$$
Equating coefficients shows that
\begin{equation}\label{eq:R2n+1 Vector Field}
X_H \ = \ H_{p_j}\partial_{q_j} - (H_{q_j}+p_jH_z)\partial_{p_j}+(p_jH_{p_j}-H)\partial_z.
\end{equation}
Or, as equations of motion ($j=1,\dots,n$),
\begin{equation}\label{eq:local on R2n+1}
\left\{\begin{array}{ccl}
\dot q_j &=& H_{p_j}\\[4pt]
\dot p_j &=& -H_{q_j}-p_jH_z \\[4pt]
\dot z &=& p_jH_{p_j}-H.
\end{array}\right.
\end{equation}

\paragraph{Not Poisson brackets}
In the more familiar symplectic setting, the Poisson brackets are defined by,
$\{H,\, f\} := X_H(f)$ -- the derivative of $f$ along the vector field $X_H$ associated to the Hamiltonian $H$. A key property (following from the skew-symmetry of the symplectic form) is that $\{g,\,f\}=-\{f,g\}$.  

In the contact setting, one can of course define a `bracket' in the same way, but it is no longer skew-symmetric. A simple calculation shows
\begin{equation}\label{eq:PoissonBrackets}
X_H(f) = \left\{H,\,f\right\}_{(p_i,q_i)} + p_i\left\{H,\,f\right\}_{(p_i,z)}-H\,f_z,
\end{equation}
where, given variables $x,y$, we write $\left\{H,\,f\right\}_{(x,y)} = H_xf_y-H_yf_x$, following the notation of \cite{Bravetti+C+T2016}.  The lack of skew-symmetry is in the final term; in particular it is a derivation of $f$ but not of $H$.  The expression for $X_H(H)$ recovers the one in \eqref{eq:derivative of H}.

\begin{remark}\label{rmk:linear}
It is perhaps a natural question to ask about linear contact vector fields on $(\RR^{2n+1},\,\eta)$ (or perhaps unnatural since $\eta$ is not linear). Using these Darboux coordinates, it is straightforward to check that the space of linear contact vector fields is only $(n^2+1)$-dimensional. With coordinates $(q_i,p_j,z)$ they have matrix
$$L= \begin{pmatrix}
A&0&0\cr 0& -A^T+aI_n&0\cr0&0&a
\end{pmatrix},$$
The Hamiltonian of this vector field is $H=p^T Aq - az $. 
where $a\in \RR$ and $A$ is any $n\times n$ real matrix. As a Lie algebra, this is isomorphic to $\mathfrak{gl}_n(\RR)\times\RR$. 
See Remark\,\ref{rmk:weights} for an extension to weighted homogeneity.
\end{remark}

%%%%%%%%%%%%%%%%%%%%%%%%%%%
%%%%%%%%%%%%%%%%%%%%%%%%%%%
\section{Equilibria}
\label{sec:equilibria}

Let $(M,\xi)$ be a contact manifold, and let $X$ be a contact vector field. 
Our study is local (in neighbourhoods of equilibria) and it will be convenient to fix a contact 1-form $\eta$, which one can always do locally.  Let $H=-\iota_X\eta$ be the associated Hamiltonian.

Equilibria occur where $X=0$. 
The definition of $X=X_H$ in \eqref{eq:Hamiltonian->VF} then yields the following conditions on the Hamiltonian function at an equilibrium point,
$$H=0,\quad \dH=\Reeb(H)\,\eta.$$
Notice that the second equation in particular implies $\dH$ is parallel to $\eta$, and if at a point $x\in M$, $\dH_x=-\tau\eta_x$ then $\Reeb(H)=\dH(\Reeb) = -\tau\eta(\Reeb)=-\tau$ so the second equation is in fact equivalent to $\dH$ being parallel to $\eta$ and hence equivalent to $\dH\restr{\xi}=0$. This shows,

\begin{proposition} \label{prop:equilibrium}
Let $(M,\xi)$ be a contact manifold and $\eta$ a (local) choice of contact form. Suppose $X\in\kX_\xi$ and  $H=-\eta(X)$ is the associated Hamiltonian function. A point $x_0\in M$ is an equilibrium point of the vector field if and only if 
$$H(x_0)=0 \quad\text{and}\quad \dH_{x_0}= -\tau\eta_{x_0},$$ 
for some $-\tau\in\RR$. In this case $\tau=-\Reeb(H)(x_0)$.
\end{proposition}

(We use $-\tau$ here to compensate for the minus sign in the definition of $H$.)
At an equilibrium point, we call the quantity $\tau=-\Reeb(H)(x_0)$ the \defn{principal coefficient} of the equilibrium (we will see below that it is an eigenvalue). 
It is not hard to show (see Proposition\,\ref{prop:invariance of reeb coeff}) that this depends only on $X$ and not on the choice of contact 1-form $\eta$. We note that for conservative contact vector fields,  the principal coefficient always vanishes. 

The proposition has a simple geometric interpretation. Namely, equilibria occur where either $H$ has a critical point on $H^{-1}(0)$ or the contact hyperplane is tangent to the zero level-set of the Hamiltonian.

\begin{proposition}\label{prop:non-singular}
Using the notation of the previous proposition, the equilibrium point $x_0$ is non-degenerate if 
\begin{enumerate}
\item the principal coefficient\, $\tau\neq0$, and
\item the bilinear form on the contact hyperplane $\xi(x_0)$ given by
$$(\uu,\vv) \ \longmapsto \ \D_\uu\left(\dH+\tau\eta\right)(\vv)$$
is non-degenerate. 
\end{enumerate}
\end{proposition}

Here we use $\D$ to denote the ordinary derivative, as distinct from the exterior derivative.  
Note that if $\alpha$ is a 1-form and $\alpha(x_0)=0$ then $D\alpha(\uu)=D_\uu\alpha$ (the derivative of $\alpha$ in the direction $\uu$) is a well-defined quantity in the cotangent space at $x_0$; that is, it is independent of any choice of coordinates.  Moreover, in any coordinates, one has 
$$\D(\dH+\tau\eta) = \D^2H + \tau \D\eta.$$
The first term of this bilinear form is the Hessian of $H$ (which does depend on coordinates, unless $H$ is singular at this point, in which case $\tau=0$).
 
\begin{definition}\label{def:amended Hessian}
We call the bilinear form $\D\left(\dH+\tau\eta\right)$ on $\xi$ at an equilibrium point the \defn{amended Hessian} of $H$, and we denote it $\Hess'$, or $\Hess'(H)$ if needed.
\end{definition}

To clarify the definition of the amended Hessian we can use local coordinates. Let $\eta=a_i\dx^i$. Then degeneracy of the bilinear form means $\exists\uu\in\xi,\;\uu\neq0$, such that
$$\forall\vv\in\xi,\quad
\left(\frac{\partial^2H}{\partial x^j\partial x^i} +\tau\frac{\partial a_i}{\partial x^j}\right)u^jv^i=0.$$
In particular, using Darboux coordinates on $\RR^{2n+1}$ about the point in question, with $\eta=\dz-p_j\dq_j$, the amended Hessian is the $2n\times 2n$ matrix, 
\begin{equation}\label{eq:amended Hessian in Darboux coords}
\Hess' = \begin{pmatrix}
H_{qq} & H_{pq}-\tau I_n \cr H_{qp} & H_{pp}
\end{pmatrix}
\end{equation}
where $H_{qq}$ is the $n\times n$ block  $H_{q_iq_j}$ etc. 
It is clear from this expression that the amended Hessian is not in general symmetric; indeed, it may have complex eigenvalues. 

For points other than the origin in Darboux coordinates, we can use the basis for $\xi$ given in \eqref{eq:basis for xi}, for which the expression for the amended Hessian becomes 
$$\Hess' = \begin{pmatrix}
H_{q_iq_j}+2p_i H_{q_jz}+H_{zz}p_ip_j & H_z\delta_{\ij}+H_{p_iq_j}+p_i H_{p_jz} \\
H_{p_jq_i}+p_i H_{p_jz} & H_{p_ip_j}
\end{pmatrix}.
$$

\begin{proof}
We prove this using the Hamiltonian; further below we see an argument using the vector field.
The equations for an equilibrium in Proposition\,\ref{prop:equilibrium} are equations of $(x,-\tau)\in M\times\RR$ (or $\RR^{2n+1}\times\RR$). 
Differentiating these equations in the direction of $\uu$ gives
$$\dH(\uu)=0,\quad \D^2H(\uu)+\tau \D_\uu\eta = -\hat\tau\eta$$
where $\D_\uu\eta$ is the derivative of $\eta$ in the $\uu$ direction, and $\hat\tau\in\RR$. 

Firstly, if $\tau=0$ then the first equation is void, and there are thus $2n+1$ equations in $2n+2$ variables, and the equations are degenerate.

However, if $\tau\neq0$ then the first equation tells us $\uu\in\xi$. For a given $\uu\in\xi$, the existence of $\hat\tau$ satisfying the second equation is equivalent to having a zero of the restriction of the linear form $(\D^2H(\uu)+\tau\D_\uu\eta)$ to $\xi$. Thus non-degeneracy is equivalent to the following bilinear form on $\xi$ being non-degenerate:
$$(\uu,\vv) \longmapsto \D^2H(\uu,\vv)+\tau \D_\uu\eta(\vv).$$
\vskip-5mm
\end{proof}

Now consider the linearization $L:T_{x_0}M\to T_{x_0}M$ of the vector field $X_H$ at an equilibrium point $x_0$.

\begin{theorem}\label{thm:linear}
Let $x_0\in M$ be an equilibrium point of $X_H$, with principal coefficient $\tau$, and let $L$ be the linear part of the vector field at $x_0$. Then, 
\begin{enumerate}
\item $L$ leaves $\xi$ invariant; we will denote the restriction to $\xi$ by $L_\xi$.
\item The linear vector field $L_\xi-\half\tau I_\xi$ on $\xi$ is Hamiltonian, where $I_\xi$ is the identity map on $\xi$,  and the symplectic structure on $\xi$ is given by the 2-form $\d\eta$;  the Hamiltonian function is given by the symmetric part of the amended Hessian. 
\end{enumerate}
\end{theorem}

The simplest proof of this statement uses the expression for $L$ in local Darboux coordinates. Since this expression will be useful later, we calculate it here. 

Differentiating the local expression \eqref{eq:local on R2n+1} for the vector field, using  Darboux coordinates $(q_j,p_j,z)$ one finds (here $j$ denotes the row and $i$ the column within each block),
$$\D(X_H) =\begin{pmatrix}
H_{p_jq_i} & H_{p_jp_i} & H_{p_jz} \\[4pt]
-H_{q_jq_i}-p_jH_{q_iz} &-H_{q_ip_j}-H_z\delta_{\ij}-p_jH_{p_iz}& -H_{q_jz}-p_jH_{zz} \\[4pt]
p_kH_{q_ip_k}-H_{q_i} & p_kH_{p_ip_k} & p_kH_{p_kz}-H_z
\end{pmatrix}.
$$ 
Evaluating this at the origin, which we assume to be an equilibrium point, gives the linear part of the vector field:
\begin{equation}\label{eq:L0 in general}
L =\begin{pmatrix}
H_{p_jq_i} & H_{p_jp_i} & H_{p_jz} \\
-H_{q_iq_j}&-H_{q_ip_j}+\tau \delta_{\ij} & -H_{q_jz} \cr
0 & 0 & \tau
\end{pmatrix}.
\end{equation}
For the record, we note that the trace of this matrix is given by
\begin{equation}\label{eq:trace}
\tr L = (n+1)\tau.
\end{equation}

It will also be useful to have the expression for the restriction $L_\xi$:
\begin{equation}\label{eq:Lxi}
L_\xi = \begin{pmatrix}
H_{p_jq_i} & H_{p_jp_i} \\
-H_{q_iq_j}&-H_{q_ip_j}+\tau \delta_{\ij} 
\end{pmatrix}.
\end{equation}

\begin{proof}
(i) The invariance of $\xi$ under $L$ follows from the zeros in the bottom row of $L$ in \eqref{eq:L0 in general} (more geometrically it follows from the fact that the flow preserves the distribution $\xi$, which also implies that $\eta(x_0)$ is an eigen-covector, or left eigenvector, of $L$; see Sec\,\ref{sec:principal coefficient} below).

(ii)  The restriction of  $L$  to $\xi$ is given in \eqref{eq:Lxi}. Hence,
$$L_\xi-\half\tau I_\xi =
\begin{pmatrix}
H_{p_jq_i} -\half\tau \delta_{\ij}& H_{p_jp_i} \\[4pt]
-H_{q_iq_j}&-H_{q_ip_j}+\half\tau \delta_{\ij} 
\end{pmatrix}.
$$
Let $J=\begin{pmatrix}0& \delta_{\ij}\cr -\delta_{\ij}&0\end{pmatrix}$ --- the matrix associated to the symplectic structure $\d\eta=\dq_i\wedge\dpp_i$ --- then
$$J\left(L_\xi-\half\tau I_\xi\right) = 
\begin{pmatrix}
H_{q_iq_j}&H_{q_ip_j}-\half\tau \delta_{\ij} \\[4pt]
H_{p_jq_i} -\half\tau \delta_{\ij}& H_{p_jp_i} 
\end{pmatrix},$$
which is the Hessian matrix at the origin of $H-\half\tau p_iq_i$, restricted to $\xi$. That is $L_\xi-\half\tau I_\xi$ is the linear vector field on $\xi$ associated to this quadratic Hamiltonian (see eg, \cite{Arnold}). 
Finally we see that this Hessian matrix is precisely the symmetric part of the amended Hessian \eqref{eq:amended Hessian in Darboux coords}. 
\end{proof}

It is well-known \cite{A+M-FoM,Arnold} that eigenvalues of a Hamiltonian (infinitesimally symplectic) matrix arise in quadruplets $\{\pm\lambda,\,\pm\bar\lambda\}$ (not necessarily all distinct). It follows from the theorem that the eigenvalues of $L_\xi-\half\tau I_\xi$ arise in these symplectic quadruplets 
and the following result is then immediate.

\begin{corollary}\label{coroll:quadruplets}
One of the eigenvalues of a contact equilibrium is equal to the principal coefficient $\tau$, while the others arise in quadruplets of the form
$$\left\{\,\half\tau\pm\lambda,\;\half\tau\pm\bar\lambda\,\right\}.$$
The eigenvalue $\tau$ corresponds to the eigen-covector $\eta$, while the others arise from the restriction to $\xi$.  
\end{corollary}

We call these \defn{contact quadruplets} of eigenvalues. Note that if $\tau\neq0$ at most 2 members of such a quadruplet may vanish or be pure imaginary. The eigenvalue $\tau$ we also call the \defn{principal eigenvalue}. 

\medskip

Recall that an equilibrium point of a vector field is non-degenerate provided the linear part has no zero eigenvalues. It follows immediately from \eqref{eq:L0 in general} that this is equivalent to,
\begin{enumerate}
\item $\tau\neq0$, and
\item in local Darboux coordinates about $x_0$, the $2n\times 2n$ matrix at $x_0$,
$$\begin{pmatrix}
H_{qq} & H_{pq}-\tau I_n \cr H_{qp} & H_{pp}
\end{pmatrix}
$$
is invertible, where we write $H_{qq}$ for the $n\times n$ matrix $(H_{q_iq_j})$ evaluated at $x_0$, etc.; this matrix is the amended Hessian (Definition\,\ref{def:amended Hessian}) in local Darboux coordinates. 
\end{enumerate}
This is equivalent to the non-degeneracy described in Proposition\,\ref{prop:non-singular}.

One sees for example that the equilibrium at the origin in $\RR^3$ for the Hamiltonian $H=z+pq$ is non-degenerate, while the one for $H=z-pq$ is degenerate (the eigenvalues can be read off the matrix $L$ in Remark\,\ref{rmk:linear}).

%%%%%%%%%%%%%%%%%%%%%%%%%%%
\subsection{Non-degenerate equilibria in  dimension 3}
Using Darboux coordinates in a neighbourhood of the origin, we consider the general Hamiltonian assuming the origin is an equilibrium point and expanded to order 2:
\begin{equation}\label{eq:Ham in R3}
H=-\tau z+ Aq^2+Bqp+ Cp^2+DDqz+EEpz+Fz^2 +O(3).
\end{equation}
Note that $\tau$ is the principal coefficient (or eigenvalue) of $H$ at the origin. 

The amended Hessian (Definition\,\ref{def:amended Hessian}) for this Hamiltonian is
$$\mathrm{Hess}':=\begin{pmatrix}
2A&B-\tau \cr B&2C
\end{pmatrix}.
$$
For non-degeneracy, we require $\tau\neq0$ and $\det\mathrm{Hess}'\neq0$. 
We remark that the eigenvalues of the amended Hessian are $(A+C)\pm\half\sqrt{4(A-C)^2+B(B-\tau)}$ which are complex if $B\tau$ is sufficiently large (positive).

The linear part of the vector field at the origin is, by \eqref{eq:L0 in general}, 
\begin{equation}\label{eq:L0 basic}
L=\begin{pmatrix}
B&2C&EE\cr
-2A & -B+\tau  & -DD\cr
0&0&\tau 
\end{pmatrix}.
\end{equation}
This has determinant $-\tau (B^2+B-\tau -4AC)$. If this is non-zero then the origin is an isolated non-degenerate equilibrium. 

The eigenvalues of $L$ are
$$\tau , \quad \frac12\left(\,\tau \pm\sqrt{(2B-\tau )^2-16AC}\,\right).$$
The equilibrium is \emph{asymptotically stable} if the real parts of the three eigenvalues all have negative real part. Thus we have (recall $\tau$ is the principal coefficient or eigenvalue),

\begin{theorem}\label{thm:stability in R3}
The origin is an asymptotically stable equilibrium if 
\begin{equation}\label{eq:asymptotic stability R3}
\tau <0\quad \text{and}\quad B^2-B\tau -4AC<0.
\end{equation}
If either of the inequalities is reversed then the equilibrium is unstable. 
\end{theorem}

%%%%%%%%%%%%%%%%%%%%%%%%%%%
\subsection{Non-degenerate equilibria in higher dimensions}

From \eqref{eq:L0 in general}, we see that $\mathrm{tr}(L)=(n+1)\tau$, so $\tau<0$ is a necessary condition for the asymptotic stability of an equilibrium point. 

As discussed in Corollary\,\ref{coroll:quadruplets}, one eigenvalue is $\tau$ and the others arise in contact quadruplets, which are of the form 
$$\left\{\,\half\tau\pm\lambda,\;\half\tau\pm\bar\lambda\,\right\}.$$

\begin{theorem}\label{thm:asymptotic stability}
Suppose an equilibrium has negative principal coefficient ($\tau<0$) and the symmetric part of the amended Hessian is positive or negative definite, then the equilibrium is asymptotically stable. 
\end{theorem}

This sufficient condition is certainly not necessary in general.

\begin{proof}
If the symmetric part of the amended Hessian is definite, then all the eigenvalues of the associated Hamiltonian system from Theorem\,\ref{thm:linear} are pure imaginary.  In this case, and with $\tau<0$, the contact quadruplets all have negative real part equal to $\tau/2$.
\end{proof}

%%%%%%%%%%%%%%%%%%%%%%%%%%%
\subsection{Hopf bifurcation}
\label{sec:Hopf}

Recall that in a dynamical system, a Hopf (or Andronov-Hopf) bifurcation occurs when a pair of eigenvalues of a non-degenerate equilibrium pass through the imaginary axis \cite{Guckenheimer-Holmes,Kuznetsov}. This gives rise to the existence of a family of periodic orbits emanating from the equilibrium point. 

Using the information above it is straightforward to construct examples of Hopf bifurcation in systems with dimension at least 5. In dimension 3 it follows from the structure of the contact quadruplets that a simple Hopf bifurcation is not possible (on the other hand a fold-Hopf bifurcation is possible --- see \S\ref{sec: Type 1 unfoldings} below). 

An explicit example in dimension 5  is to let
$$H_\lambda =z + p_1q_2 - q_1p_2 + 2\lambda q_1p_1+O(3).$$
The origin is an equilibrium point with principal coefficient $-1$ (for all $\lambda$); the corresponding principal coefficient is $-1$, while the other eigenvalues are
$$\lambda \pm\sqrt{-1+\lambda^2} ,\quad  -1 -\lambda \pm \sqrt{-1+\lambda^2}.$$
When $\lambda=0$ these form the contact quadruplet $\{\pm i,\; -1\pm i\}$. As $\lambda$ varies, the first two cross the imaginary axis with non-zero velocity (their real part is equal to $\lambda$), as shown in Figure\,\ref{fig:Hopf}.  Without the $O(3)$ terms, this system is linear (see Remark\,\ref{rmk:linear}) and this would give rise to a `vertical' Hopf bifurcation, meaning that the periodic orbits all occur for $\lambda=0$ (in fact in the $q_1q_2$ plane). The addition of suitable higher order terms  would make it a sub- or super-critical Hopf bifurcation. Note that for $\lambda<0$ (small) the origin is asymptotically stable, while for $\lambda>0$ it is unstable.  

\begin{figure}[t] % Hopf 
\psset{unit=0.5}
\centering
\begin{pspicture}(-1,-3)(2,2)
\psline[linecolor=Grey]{->}(-2,0)(2,0) 
\psline[linecolor=Grey]{->}(0,-2)(0,2) 
\psdot[linecolor=Grey](-1,0)
\psdots(-0.7,0.953)(-0.7,-0.953)(-0.3, 0.953)(-0.3, -0.953)
\rput(0,-2.4){\footnotesize$\lambda<0$}
\end{pspicture}
\qquad\qquad
\begin{pspicture}(-1,-3)(2,2)
\psline[linecolor=Grey]{->}(-2,0)(2,0) 
\psline[linecolor=Grey]{->}(0,-2)(0,2) 
\psdot[linecolor=Grey](-1,0)
\psdots(-1,1)(-1,-1)(0,1)(0,-1)
\rput(0,-2.4){\footnotesize$\lambda=0$}
\end{pspicture}
\qquad\qquad
\begin{pspicture}(-2,-3)(2,2)
\psline[linecolor=Grey]{->}(-2,0)(2,0)  
\psline[linecolor=Grey]{->}(0,-2)(0,2)  
\psdot[linecolor=Grey](-1,0)
\psdots(-1.3,0.9)(-1.3,-0.9)(0.3, 0.953)(0.3, -0.953)
\rput(0,-2.4){\footnotesize$\lambda>0$}
\end{pspicture}
\caption{Contact quadruplet exhibiting a Hopf bifurcation in dimension 5}
\centerline{(the grey dot represents the principal coefficient) --- see \S\,\ref{sec:Hopf}.}
\label{fig:Hopf}
\end{figure}

%%%%%%%%%%%%%%%%%%%%%%%%%%%
\subsection{Degenerate equilibria}
\label{sec:degeneracies}
It follows from Proposition \ref{prop:non-singular} that there are two distinct ways in which an equilibrium can be degenerate. We call them Type I and Type II degeneracies as follows.
 
\begin{definition}\label{def:Types} A degenerate equilibrium $x_0$ with simple zero eigenvalue, of a contact Hamiltonian system is of
\begin{itemize}
\item \defn{Type I}\ \ if the principal coefficient vanishes (in this case $H$ has a critical point at $x_0$);
\item \defn{Type II}\ \ if the amended Hessian matrix is degenerate.
\end{itemize}
\end{definition}

We will see that these are generically of codimension 1.  Higher codimension degeneracies can occur that combine the two types, but we do not study these in this paper.

There follow below two parallel sections, one on Type I singularities and the other on Type II singularities. The analysis of the first type is the more straightforward, because the condition for a fold singularity only depends on the 2-jet of the Hamiltonian, while for Type II it depends on its 3-jet. In each section, we begin with a general discussion of the singularities in $\RR^{2n+1}$ and then follow it with a section on the 3-dimensional cases.  Before we proceed with that analysis, we address the question of the dependence of the principal coefficient and amended Hessian on the choice of contact form. 

\subsection{Dependence on \texorpdfstring{$\eta$}{eta}}    Recall that, given a contact manifold $(M,\xi)$ and a contact vector field $X$, the Hamiltonian itself depends on the choice of contact form $\eta$.  We show directly that the principal coefficient of $X$ at an equilibrium point is independent of the choice of contact form, and the amended Hessian is well-defined up to scalar multiple (although the first part also follows from the fact that the principal coefficient is $-2$ times the principal coefficient of the vector field).  

\begin{proposition}\label{prop:invariance of reeb coeff}
Let $X$ be a contact vector field on the contact manifold $(M,\xi)$, and let $x_0\in M$ be an equilibrium point. 
\begin{enumerate}
\item The principal coefficient $\tau$ of\/ $X$ at $x_0$ is independent of the choice of contact form. 
\item The amended Hessian is well-defined up to scalar multiple.  More precisely, if $\eta_1,\eta_2$ are two 1-forms representing $\xi$, so that $\eta_2=f\eta_1$ for some non-zero function $f$, then the bilinear forms on $\xi$ at the equilibrium point $x_0$ satisfy\/ $\Hess_2'= f(x_0)\, \Hess_1'$. 
\end{enumerate}
\end{proposition}

\begin{proof}
(i)  Since $\tau$ is an eigenvalue of the vector field, it does not depend on any choice arising from the contact form $\eta$.  It also follows from the fact that $\tr(L)=(n+1)\tau$, see \eqref{eq:trace}. 

(ii) For $\uu,\vv\in\xi$, the amended Hessians are defined by
$$\Hess_j'(\uu,\vv) = \D_\uu(\dH_j +\tau\eta_j)(\vv).$$
Now, $\eta_2=f\eta_1$  implies $H_2 = f H_1$, and hence
$$\dH_2 = f\,\dH_1 + H_1\df$$
and then
$$\dH_2+\tau\eta_2 = f(\dH_1+\tau\eta_1)+H_1\df.$$
Thus  (in any coordinate system), 
$$\D\left(dH_2+\tau\eta_2\right) = f\,\D(\dH_1+\tau\eta_1) + 
\df\otimes(\dH_1+\tau\eta_1) +\dH_1\otimes\df + H_1\D^2f.$$
Then, at the equliibrium point $x_0$ and restricting to $\xi$, all but the first term vanishes, showing that 
$$\Hess_2'=f(x_0)\,\Hess_1'$$
as required. 

\end{proof}

\def\wt{\mathop\mathrm{wt}}
\begin{remark}\label{rmk:weights} 
Perhaps more natural than Remark\,\ref{rmk:linear} is to assign weights to the Darboux coordinates: 
$$\wt(q_i)=\wt(p_i)=1,\;\wt(z)=2.$$
Then $\eta$ is homogenous of degree 2. And given a Hamiltonian function which is homogeneous of degree $d$ then the vector field has degrees $d-1$ in the first $2n$ components and $d$ in the last component, meaning that the vector field itself is homogeneous of degree $d-2$.  
From \eqref{eq:PoissonBrackets} one sees that if $H$ and $f$ are weighted homogeneous, with $\deg(H)=d$, $\deg(f)=r$, then $X_H(f)$ has degree $r+d-2$.  For example, let $H$ be the general Hamiltonian of weighted degree 2 in $\RR^3$,
$$H=-\tau z + Aq^2+Bqp+Cp^2$$
($A,B,C,\tau\in\RR$), then the corresponding vector field is,
$$X_H= \begin{pmatrix}
2 B q +2 C p  \\
 \left(-2 B +\tau \right) p -2 A q  \\
 -A \,q^{2}+C \,p^{2}+\tau  z 
\end{pmatrix}.
$$
which is of weighted degree 0.
\end{remark}

\subsection{Principal coefficients of contact diffeomorphisms}
\label{sec:principal coefficient}
Here we remark on a geometric view of the principal coefficient. 

Given a contact manifold $(M,\xi)$, denote by $\xi^\circ \subset T^*M$ the line bundle of linear forms vanishing on $\xi$ (that is, $\xi^\circ$ is the annihilator of $\xi$). Now any diffeomorphism $\Phi$ of $M$ preserving $\xi$ will also preserve $\xi^\circ$.  

\begin{definition}\label{def:principal coefficient FP}
Suppose $x_0\in M$ is a fixed point of such a contact diffeomorphism. Then it (or rather its cotangent lift) maps $\xi^\circ(x_0)$ to itself, acting by scalar multiplication. We call the corresponding scalar the \defn{principal coefficient} of the diffeomorphism at the fixed point. 
\end{definition}

Recall that if $\Phi$ is a diffeomorphism then the cotangent lift $\Phi^*$ is given by 
$$\left<\Phi^*(\alpha_y),\,v_x\right> := \left<\alpha_y,\,\d\Phi_x(v_x)\right>,$$
where $y=\Phi(x)$, $\alpha_y\in T^*_yM$ and $v_x\in T_xM$. 

If, in a neighbourhood of a fixed point, we chose a contact 1-form $\eta$, then the contact diffeomorphism maps $\eta$ to another contact 1-form, which is of the form $f\,\eta$, where $f$ is a non-vanishing function; that is $\Phi^*\eta=f\,\eta$.  The principal coefficient of the diffeomorphism $\Phi$ at a fixed point $x_0$ is then just $f(x_0)$. This value is clearly independent of the choice of contact form.

Now suppose $X\in\kX_\xi$ is a contact vector field and $\Phi_t$ its flow. Let $x_0$ be an equilibrium point of $X$. Choosing an arbitrary contact form $\eta$ in a neighbourhood of $x_0$, let $\Phi^*_t\eta=f_t\,\eta$. Let $\tau$ be the principal coefficient of $X$. Then we have 
$f_t(x_0) = \exp(t\tau)$, and 
$$\tau  =  \frac{\d}{\dt}\, f_t(x_0)\restr{t=0}.$$

%%%%%%%%%%%%%%%%%%%%%%%%%%%
%%%%%%%%%%%%%%%%%%%%%%%%%%%
\section{Degeneracy of Type I}
A type I  degeneracy of an equilibrium is one where the principal coefficient  vanishes; in other words it arises at a singular point of the zero-level of $H$.  We consider now the conditions for this to be a simple degeneracy (i.e., a fold).  

First we assume the zero eigenvalue of the linear part of the vector field at the equilibrium is simple, and then we ask when the singularity is of fold type.  

\subsection{Type I fold singularity}
\label{sec:Type 1 fold singularity}
With $\tau=0$, the linear part $L$ takes the form (see \eqref{eq:L0 in general})
\begin{equation}\label{eq:L0 with tau=0}
L =\begin{pmatrix}
H_{pq} & H_{pp} & H_{pz} \\
-H_{qq}&-H_{qp} & -H_{qz} \cr
0 & 0 & 0
\end{pmatrix}.
\end{equation}
Here $H_{qq}=\left(H_{q_iq_j}\right)$ etc.

This clearly has corank at least 1, and as stated above, we begin by requiring that 0 is a simple eigenvalue, which in particular means the matrix has corank 1. More precisely it requires that the top left $2n\times 2n$ block $L_\xi$ be non-degenerate (equivalently, the Hessian of the Hamiltonian on $\xi$ be non-degenerate).

\begin{theorem}\label{thm:Delta2}
Suppose a contact dynamical system with Hamiltonian $H$ has a degenerate equilibrium with vanishing principal coefficient. The singularity is a fold if,
$$\Delta_2:=H_{zq}\aa_q+H_{zp}\aa_p+H_{zz}\neq0$$
where \ $\aa_q,\aa_p\in\RR^n$ satisfy
\begin{equation}\label{eq:x and y}
\left\{\begin{array}{rcl}
H_{qq}\aa_q+H_{qp}\aa_p+H_{qz}&=&0,\\[4pt]
H_{pq}\aa_q+H_{pp}\aa_p+H_{pz}&=&0.
\end{array}\right.
\end{equation}
Furthermore, in this case, the family $H_\lambda=H-\lambda$ is a versal unfolding of the singularity (a saddle-node bifurcation). 
\end{theorem}

Explicitly, the $i^{\textrm{th}}$ component of the first equation in \eqref{eq:x and y} is
$$\sum_j\left(H_{q_iq_j}(\aa_q)_j+H_{q_ip_j}(\aa_p)_j\right)+H_{q_iz}=0.$$
The components of the second equation and the expression for $\Delta_2$  are similar. 

\begin{proof} 
We apply Lemma\,\ref{lemma:fold} from the appendix.  To do so we need non-zero vectors $\vv\in\mathrm{coker} L$ and $\aa\in\ker L$ to check whether $\vv \, \D^2(X_H)\,\aa^2\neq0$.  

Now $\vv=(0,\,0,\dots,0,1)$ is clearly a non-zero element of the cokernel (ie, $\vv\,L=0$). 

To find $\aa\in\ker L$ we know $\aa\not\in\xi$ so we can choose it of the form $\aa=(\aa_q,\,\aa_p,\,1)^T$ with $\aa_q,\aa_p\in\RR^n$. Then \eqref{eq:x and y} is precisely the condition that $\aa\in\ker L$. 

By Lemma\,\ref{lemma:fold}, we require $\vv \,\D^2(X_H)\,\aa^2\neq0$. Now $\vv\, \D^2(X_H)$ is the Hessian matrix of the final component $\vv\, X_H$ of $X_H$, so let $f=\vv\, X_H = p_jH_{p_j}-H$ (a Legendre transform of $H$).  Then, at the origin, one finds
$$\D^2f = \begin{pmatrix}
-H_{qq} & 0 & -H_{zq} \\
0 & H_{pp} & 0\\
-H_{qz} & 0 & -H_{zz}
\end{pmatrix}.
$$
In order that $\vv \,\D^2(X_H)\aa^2\neq0$, we require
$$\aa^T\,\left[\D^2f\right]\,\aa\neq0.$$
Expanding that in terms of $\aa_q$ and $\aa_p$, and simplifying using \eqref{eq:x and y} gives the result.
\end{proof}

%%%%%%%%%%%%%%%%%%%%%%%%%%%
\subsection{Fold singularity in \texorpdfstring{$\RR^3$}{R3}}
We translate the condition of the theorem above to conditions on the coefficients in the Taylor series for $H$. In this case (Type I) the theorem above shows we only need the Taylor series to order 2 (at an equilibrium point).  

With vanishing principal coefficient, the lowest order terms of the Hamiltonian at an equilibrium point (in Darboux coordinates) are quadratic:
\begin{equation}\label{eq:Type I Hamiltonian}
H_0 =  Aq^2+Bqp+Cp^2 + Dqz + Epz+Fz^2 + O(3).
\end{equation}
Recall that the linearization at the origin is (with $\tau=0$) 
$$L = \begin{pmatrix}
B & 2C & E\\
-2A & -B & -D\\
0&0&0
\end{pmatrix}.
$$

Consider the two quantities,
\begin{equation}\label{eq:NDGC TypeI}
\left\{\begin{array}{ll}
\Delta_1:=&B^2-4AC\\
\Delta_2:=&4(B^2-4AC)F + AE^2 - BDE + CD^2.\\
\end{array}\right.
\end{equation}

If $\Delta_1\neq0$ then 0 is a simple eigenvalue of $L$. We assume this from now on.  

\begin{corollary}\label{coroll:Type1 NDC} Consider a contact vector field on $\RR^3$ with a singularity of type I at the origin, and hence with Hamiltonian as above, with $\Delta_1\neq0$.\\
(i) The vector field has a fold singularity if and only if $\Delta_2\neq0$. \\
(ii) In this case, the family $H_\lambda=H_0-\lambda$ gives a versal unfolding of the vector field, resulting in a saddle-node bifurcation of equilibria.
\end{corollary}

This is a particular case of Theorem\,\ref{thm:Delta2} above, but note that the $\Delta_2$ here is $(B^2-4AC)$ times the $\Delta_2$ defined in the theorem.

%%%%%%%%%%%%%%%%%%%%%%%%%%%
\subsection{Type I saddle-node bifurcations} \label{sec: Type 1 unfoldings}
Consider the 1-parameter family of Hamiltonian functions
$$H_\lambda= -\lambda + Aq^2+Bpq+Cp^2 +Dqz+Epz+Fz^2+O(3);$$ 
here $\lambda$ is the parameter, and $A,Q,\dots,F$ are fixed and satisfy $\Delta_1\neq0,\;\Delta_2\neq0$. At $\lambda=0$ this has a degenerate equilibrium of type I at the origin. The linearization $L_0$ at the bifurcation point is given in \eqref{eq:L0 basic} but with $\tau=0$, and the (amended) Hessian is
$$\Hess' = \begin{pmatrix}
2A&B\cr B&2C
\end{pmatrix}
$$
which we are assuming is non-degenerate ($\Delta_1\neq0$).

To simplify calculations, let us consider the cases where $D=E=0$ and $F=1$.  Similar results hold more generally. Then
$$H_\lambda = -\lambda + Aq^2+Bqp+Cp^2+z^2.$$
The vector field is
$$X_\lambda= \begin{pmatrix}
B q +2 C p  \\
 -2 A q -B p -2 p z  \\
-A \,q^{2}+C \,p^{2}-z^{2}+\lambda
\end{pmatrix}
$$
There are no equilibria (near 0) for $\lambda<0$ and two for $\lambda>0$:
\begin{description}
\item $(q,p,z)=(0,0,\sqrt{\lambda})$: the principal coefficient is $\tau=-2\sqrt{\lambda}$, and the eigenvalues of the linear part at the equilibrium point are 
$$+\tau, \; +\half\tau\pm\half\sqrt{(B+\tau)^2-4AC}$$ 
(with $\tau=-2\sqrt{\lambda}<0$). 
\item $(q,p,z)=(0,0,-\sqrt{\lambda})$: the principal coefficient is $\tau=2\sqrt{\lambda}$, and the eigenvalues of the linear part are as before, $\tau, \; \half\tau\pm\half\sqrt{(B-\tau)^2-4AC}$  (with $\tau>0$). This is therefore an unstable equilibrium. 
\end{description}

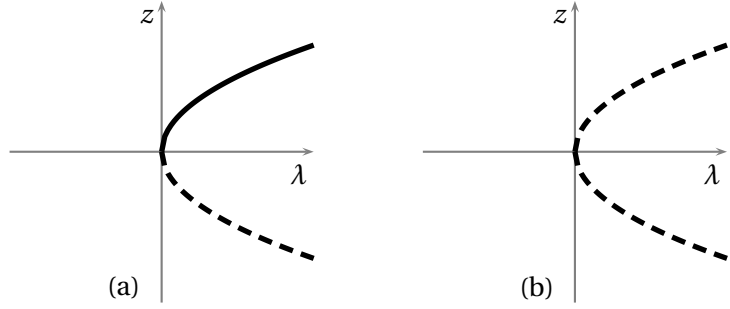
\begin{figure}[t] % Saddle-node
\centering
\begin{pspicture}(-1,-2)(2,2)
\psline[linecolor=Grey]{->}(-2,0)(2,0) \rput(-0.2,1.8){$z$}
\psline[linecolor=Grey]{->}(0,-2)(0,2) \rput(1.8,-0.3){$\lambda$}
\psplot[linewidth=2pt]{0}{2}{x sqrt}
\psplot[linewidth=2pt,linestyle=dashed]{0}{2}{x sqrt -1 mul}
\rput(-0.5,-1.8){(a)}
\end{pspicture}
\qquad\qquad
\begin{pspicture}(-2,-2)(2,2)
\psline[linecolor=Grey]{->}(-2,0)(2,0)  \rput(-0.2,1.8){$z$}
\psline[linecolor=Grey]{->}(0,-2)(0,2) \rput(1.8,-0.3){$\lambda$}
\psplot[linewidth=2pt,linestyle=dashed]{0}{2}{x sqrt}
\psplot[linewidth=2pt,linestyle=dashed]{0}{2}{x sqrt -1 mul}
\rput(-0.5,-1.8){(b)}
\end{pspicture}
\caption{Type I bifurcations: (a) the elliptic case, (b) the hyperbolic case. A solid curve represents an asymptotically stable equilibrium, and a dashed curve represents an unstable equilibrium.}
\label{fig:saddle-node}
\end{figure}

This is a saddle-node bifurcation, but there are two cases to consider:
\begin{enumerate}
\item $B^2<4AC$ (`elliptic'): for $\lambda=0$ the non-zero eigenvalues are pure imaginary. As $\lambda$ is varied, their real parts become $-2\sqrt{\lambda}$ and $-\sqrt{\lambda}$ along the `top' branch ($z>0$), and $2\sqrt{\lambda}$ and $\sqrt{\lambda}$ along the `bottom' branch  ($z<0$). The top branch of equilibria are therefore asymptotically stable, while the equilibria on the bottom branch are unstable. Note that if we changed to $F=-1$ then the equilibria would occur for $\lambda<0$, but otherwise the analysis would be unchanged. 

Moreover, passing through $0$ along the curve of equilibria, two eigenvalues cross the imaginary axis showing this is a fold-Hopf bifurcation which is normally a codimension 2 bifurcation, see Guckenheimer and Holmes \cite[Sec 7.4]{Guckenheimer-Holmes} and Kuznetsov \cite[Sec 8.5]{Kuznetsov}, but here this is exhibited as a codimension 1 phenomenon. Which dynamical phenomena are associated to this bifurcation needs further consideration --- presumably different values of the coefficients will lead to different paths through the generic co\-dimension-2 fold-Hopf bifurcation described in \cite{Guckenheimer-Holmes, Kuznetsov}. 

\item $B^2>4AC$ (`hyperbolic'): for $\lambda=0$ the non-zero eigenvalues are real, one positive, one negative.  As $\lambda$ is varied, their signs don't change and the bifurcating equilibria are therefore both unstable, and of the two equilibria one will have 1 negative and 2 positive eigenvalues while the other has  1 positive and 2 negative eigenvalues.  
\end{enumerate}

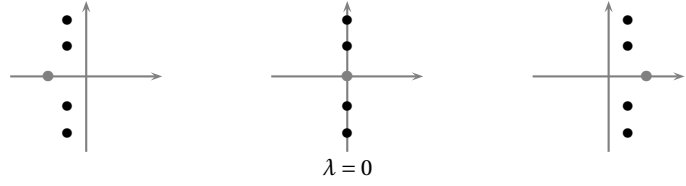
\begin{figure}[t] % double Hopf 
\psset{unit=0.5}
\centering
\begin{pspicture}(-2,-3)(2,2)
\psline[linecolor=Grey]{->}(-2,0)(2,0) 
\psline[linecolor=Grey]{->}(0,-2)(0,2) 
\psdot[linecolor=Grey,dotsize=4pt](-1,0)
\psdots(-0.5,0.8)(-0.5,-0.8)(-0.5, 1.5)(-0.5, -1.5)
\end{pspicture}
\qquad\qquad
\begin{pspicture}(-2,-3)(2,2)
\psline[linecolor=Grey]{->}(-2,0)(2,0) 
\psline[linecolor=Grey]{->}(0,-2)(0,2) 
\psdot[linecolor=Grey,dotsize=4pt](0,0)
\psdots(0,0.8)(0,-0.8)(0, 1.5)(0, -1.5)
\rput(0,-2.4){\footnotesize$\lambda=0$}
\end{pspicture}
\qquad\qquad
\begin{pspicture}(-2,-3)(2,2)
\psline[linecolor=Grey]{->}(-2,0)(2,0)  
\psline[linecolor=Grey]{->}(0,-2)(0,2)  
\psdot[linecolor=Grey,dotsize=4pt](1,0)
\psdots(0.5,0.8)(0.5,-0.8)(0.5, 1.5)(0.5, -1.5)
\end{pspicture}
\caption{Contact quadruplet for a Type I degeneracy exhibiting a fold+double-Hopf bifurcation in dimension 5. The grey dot represents the principal coefficient eigenvalue and is equal to twice the real part of the other eigenvalues. If we ignore the extreme dots, we would have a 3-dimensional fold-Hopf bifurcation. 
}
\label{fig:double-Hopf}
\end{figure}

\paragraph{In 5 and more dimensions}
A similar analysis in dimension 5 or more allows for an elliptic case, where the Hessian of the Hamiltonian on the contact plane is positive or negative definite. In this case each `quadruplet' of eigenvalues of $L_0$ will be pure imaginary, and as one moves along the saddle-node curve the eigenvalues will generically  move across the imaginary axis. This would be a fold-multi-Hopf bifurcation, which has not been analyzed. It would usually be a codimension 3 phenomenon in $\RR^5$ (or codimension $n+1$ in $\RR^{2n+1}$), but in this contact setting arises as codimension 1. See Figure\,\ref{fig:double-Hopf}.  An added complication could arise if there are any resonances between the imaginary eigenvalues when $\lambda=0$.

\begin{GR}\label{rmk:elliptic deformation}
Suppose that the Hessian $\D^2H_0$ is positive definite at the Type I equilibrium $x_0$. Then (at least in a neighbourhood of $x_0$), the zero-set of the Hamiltonian $H_0$ is just the one point, and the positive level sets of  $H_0$ are (topologically) spheres.  It follows that the zero level-set of $H_\lambda=H_0-\lambda$ is one of those spheres when $\lambda>0$ is fixed (and small).  
As already remarked, equilibria occur at points where the contact hyperplane is tangent to the sphere $H_\lambda^{-1}(0)$.  If there were no equilibria on the sphere then there would be a nowhere vanishing vector field on the sphere, which is impossible for topological reasons since the sphere has even dimension and its Euler characteristic is 2.  Therefore there must be equilibria for each $\lambda>0$ (sufficiently small), as we have seen by direct calculation in the 3-dimensional case.  On the other hand, when $\lambda<0$ the zero level-set is empty, at least near $x_0$ so there are no equilibria.  

A similar argument applies if the Hessian is negative definite,  changing the sign of $\lambda$. 
\end{GR}

\begin{remark}\label{rmk:sandwhich}
Bravetti et al.\  \cite{Bravetti-deLeon+} consider the dynamics on $S=H^{-1}(0)$ under the assumption that $\Reeb(H)\neq0$ along that hypersurface $S$. Since $\Reeb(H)=0$ at a Type I degeneracy, it would be interesting to understand how this degeneracy and its associated saddle-node bifurcation influences their findings. 
\end{remark}

%%%%%%%%%%%%%%%%%%%%%%%%%%%
%%%%%%%%%%%%%%%%%%%%%%%%%%%
\section{Degeneracy of Type II}
Here we consider degenerate equilibria with non-zero principal coefficent. In this case the degeneracy is in the restriction $L_\xi$ of $L$ to $\xi$. 

\begin{figure}[t] % Type II 
\psset{unit=0.6}
\centering
{\begin{pspicture}(-2,-3)(2,2)
\psline[linecolor=Grey]{->}(-2,0)(2,0) 
\psline[linecolor=Grey]{->}(0,-2)(0,2) 
\psdot[linecolor=Grey,dotsize=4pt](-1.4,0)
\psdots(-0.5,0)(-0.9,0)
\psline{|->}(-0.9,0.4)(-1.6,0.4)
\psline{|->}(-0.5,0.4)(0.2,0.4)
\rput(0,-2.4){\footnotesize branch 1}
\rput(0,-3){\footnotesize (stable)}
\end{pspicture}}
\qquad\qquad
\begin{pspicture}(-2,-3)(2,2)
\psline[linecolor=Grey]{->}(-2,0)(2,0) 
\psline[linecolor=Grey]{->}(0,-2)(0,2) 
\psdot[linecolor=Grey,dotsize=4pt](-1.4,0)
\pscircle(-1.4,0){0.25} \psdot(0,0)
\rput(0,-2.6){\footnotesize$\lambda=0$}
\end{pspicture}
\qquad\qquad
\begin{pspicture}(-2,-3)(2,2)
\psline[linecolor=Grey]{->}(-2,0)(2,0)  
\psline[linecolor=Grey]{->}(0,-2)(0,2)  
\psdot[linecolor=Grey,dotsize=4pt](-1,0)
\psdots(-1.8,0)(0.8,0)
\rput(0,-2.4){\footnotesize branch 2}
\rput(0,-3){\footnotesize (unstable)}
\end{pspicture}

\begin{pspicture}(-2,-3.2)(2,3)
\psline[linecolor=Grey]{->}(-2,0)(2,0) 
\psline[linecolor=Grey]{->}(0,-2)(0,2) 
\psdot[linecolor=Grey,dotsize=4pt](-1.4,0)
\psdots(-0.5,0)(-0.9,0)
\psdots(-0.7,1)(-0.7,-1)
\rput(0,-2.4){\footnotesize branch 1}
\rput(0,-3){\footnotesize (stable)}
\end{pspicture}
\qquad\qquad
\begin{pspicture}(-2,-3.2)(2,3)
\psline[linecolor=Grey]{->}(-2,0)(2,0) 
\psline[linecolor=Grey]{->}(0,-2)(0,2) 
\psdot[linecolor=Grey,dotsize=4pt](-1.4,0)
\pscircle(-1.4,0){0.25} \psdot(0,0)
\psdots(-0.7,1)(-0.7,-1)
\rput(0,-2.6){\footnotesize$\lambda=0$}
\end{pspicture}
\qquad\qquad
\begin{pspicture}(-2,-3.2)(2,3)
\psline[linecolor=Grey]{->}(-2,0)(2,0)  
\psline[linecolor=Grey]{->}(0,-2)(0,2)  
\psdot[linecolor=Grey,dotsize=4pt](-1.4,0)
\psdots(-1.8,0)(0.8,0)
\psdots(-0.7,1)(-0.7,-1)
\rput(0,-2.4){\footnotesize branch 2}
\rput(0,-3){\footnotesize (unstable)}
\end{pspicture}
\caption{Typical `motion' of eigenvalues through a saddle-node bifurcation of Type II in 3 dimensions (above) and an example in 5 dimensions (below), both with $\tau<0$. The grey dot is the principal coefficient. Reflecting in the imaginary axis would show the typical motion for $\tau>0$ where all equilibria would be unstable.}
\label{fig:Type 2 eigenvalues}
\end{figure}
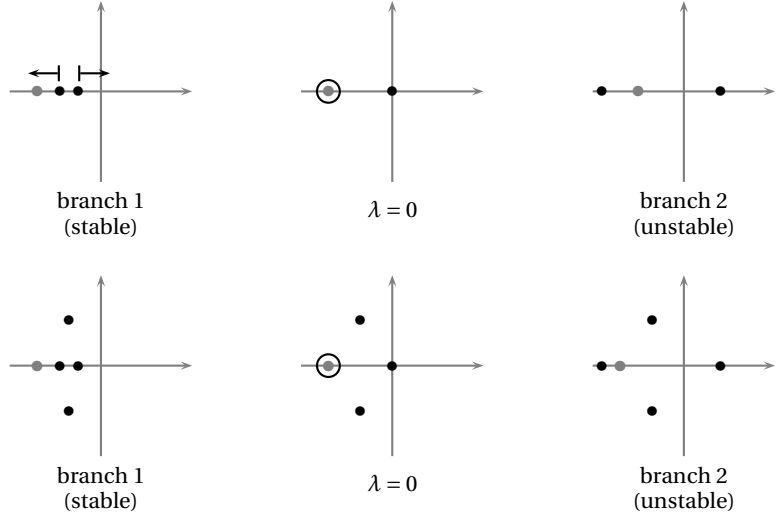

At the level of eigenvalues, we are assuming the principal eigenvalue $\tau\neq0$ and there is a (simple) zero eigenvalue. This means that one of the contact quadruplets $\left\{\,\half\tau\pm\lambda,\;\half\tau\pm\bar\lambda\,\right\}$ contains zero. This implies $\lambda=\pm\tau/2$. Then the quadruplet is simply $\{\tau,0\}$. Therefore at a degenerate equilibrium of Type II, $\tau$ is a double eigenvalue. See Figure\,\ref{fig:Type 2 eigenvalues}. However, it is not possible for the double eigenvalue to become a complex conjugate pair, as the principal eigenvalue always remains real.

\subsection{Type II fold singularity}
\label{sec:Type 2 fold singularity}

In Darboux coordinates, we saw in \eqref{eq:L0 in general} that the linear approximation at the origin is
$$L=\begin{pmatrix}
L_\xi&\rho\cr 0 &\tau
\end{pmatrix}$$
where $L_\xi$ and $\rho$ are the $2n\times 2n$ matrix and $2n$-vector,
$$L_\xi=\begin{pmatrix}
H_{pq}&H_{pp}\cr
-H_{qq}&-H_{qp}+\tau I_n  \end{pmatrix}, \qquad
\rho = \begin{pmatrix}
H_{pz}\cr -H_{qz}
\end{pmatrix}.
$$
Since $\tau\neq0$, for a degenerate equilibrium we need that $\det L_\xi=0$, and for zero to be a simple eigenvalue we require $L_\xi$ to have rank $2n-1$.

We will apply Lemma\,\ref{lemma:fold} to find conditions that ensure this has a fold singularity. In order to do this we need elements of the kernel and cokernel of $L$. 
Let $(\aa_q,\aa_p)\in\RR^{n}\times\RR^n$ and $(\vv_p,-\vv_q)\in(\RR^{n}\times\RR^n)^*$ be such that  $$(\vv_p,-\vv_q) L_\xi=0\quad\text{and}\quad L_\xi\begin{pmatrix}
\aa_q\cr\aa_p\end{pmatrix}=0.$$ 
Then $\vv L=0$ and $L\aa=0$ where
$$\vv=(\vv_p,-\vv_q,\zeta) \quad \text{and}\quad \aa = \begin{pmatrix}
\aa_q\cr \aa_p \cr0\end{pmatrix},$$
and $\zeta=-\frac1{\tau}(\vv_pH_{pz}+\vv_qH_{qz})$. 

To ensure this is a fold, rather than a more degenerate singularity, Lemma\,\ref{lemma:fold} says we need 
\begin{equation}\label{eq:NDC Type 2}
\vv\,\D^2(X_H)\,\aa^2 \neq0.
\end{equation}
This condition is equivalent to $\D_\aa^2(\vv X_H)\neq0$ (recall $\vv$ is a fixed covector). 
Unlike the Type I case, this depends on the 3-jet of the Hamiltonian at the equilibrium point. 
Written out in terms of partial derivatives, we require
$$\D_\aa^2\left(\vv_qH_q-(\vv_qp)H_z) + \vv_pH_p+\zeta(pH_p-H)\right)\neq0.
$$
which expands to (after evaluating at the origin)
\begin{equation}\label{eq:Type II NDC}
\begin{array}{llcl}
&\left(\vv_qH_{qqq} + \vv_pH_{pqq}-\zeta\,H_{qq} \right)\aa_q^2 \\[4pt]
&\qquad + 2\left(\vv_qH_{qqp}+\vv_pH_{pqp} \right)\aa_q\aa_p \\[4pt]
&\qquad\quad + \left(\vv_qH_{qpp}  +\vv_pH_{ppp}\right)\aa_p^2 \\[4pt]
&\qquad\qquad - (\vv_q\aa_p)\,\left(H_{zq}\aa_q  +H_{zp}\aa_p\right) &\neq & 0.
\end{array} 
\end{equation}
The notation should be self-explanatory. For example, with summation over repeated indices understood ($i,j,k=1,\dots,n$),
$$\vv_pH_{pqp}\aa_q\aa_p = (\vv_p)_k\,\left(\frac{\partial^3H}{\partial p_k\,\partial q_i\,\partial p_j}\right) \,(\aa_q)_i\,(\aa_p)_j.$$
This proves the first part of the following theorem.

\begin{theorem}\label{thm:Type 2 fold}
Consider an equilibrium point in $\RR^{2n+1}$ with a Type II degeneracy; that is, $\tau\neq0$ and $\mathrm{rank}(L_\xi)=2n-1$. Then, using the notation introduced above,
 \begin{enumerate}
 \item the vector field has a fold singularity provided condition \eqref{eq:Type II NDC} holds, and
 \item in this case the family $H_\lambda=H-\lambda(\alpha q+\beta p+\gamma)$, with $\alpha,\beta\in\left(\RR^n\right)^*$ and $\gamma\in\RR$, gives a versal unfolding of the singularity of the vector field, resulting in a saddle-node bifurcation of equilibria, provided
 $$(\beta,-\alpha,-\gamma)^T\not\in\Image(L).$$
 \end{enumerate}
 \end{theorem}

\begin{proof}
Part (i) is already proved by the calculation above.  

\noindent(ii) For the given Hamiltonian $H_\lambda$,
$$X_{\lambda} \ = \ X_0 + \lambda \begin{pmatrix}
\beta\cr -\alpha\cr -\alpha q -\gamma \end{pmatrix},$$
where $X_\lambda$ is the vector field associated to $H_\lambda$. 
It follows that the velocity of the deformation satisfies  $\dot{X}_\lambda(0) = (\beta,\;-\alpha,\;-\gamma)^T$ and hence the statement follows from Lemma \ref{lemma:fold}(ii) in the appendix.
\end{proof}

For part (ii), if $(H_{zq},H_{zp})\neq(0,0)$ (i.e., $\rho\neq0$) we find $H_\lambda=H-\lambda$ is a versal unfolding of the fold singularity (similar to the Type I case), whereas if $\rho=0$ it is not versal.

%%%%%%%%%%%%%%%%%%%%%%%%%%%
\subsection{Fold singularity in \texorpdfstring{$\RR^3$}{R3}}
Suppose the origin in $\RR^3$ is a degenerate equilibrium point of type II. In this case we can write the 3-jet of the Hamiltonian at the origin as
\begin{equation}\label{eq:Type2 3-jet}
\begin{split}
H \ = \ &-\tau z+ Aq^2+Bqp+ Cp^2+Dqz+Epz+Fz^2 + {}\\
&\quad + \sum_{i\leq j\leq k} P_{i,j,k}\,x_i\,x_j\,x_k \ + \ O(4),
\end{split}
\end{equation}
where $\tau\neq0$ and in the cubic terms, $x_1=q,\,x_2=p,\;x_3=z$.

Define the following 3 polynomials in the coefficients of the 3-jet of $H$ at the origin (an equilibrium point), where we write $B_1=B-\tau$:
\begin{equation}\label{eq:Type2 NDC}
\left\{\begin{array}{rcl}
	h_0 &=& B\,B_1-4AC,\\[8pt] 
	h_1&=& B^{2}(\,3 B E -6 C D - E\tau\,) \\[4pt]
	&& \quad {}+ 24 \,B \,C^{2}\,  P_{1,1,1} - 4\,B \,C\, \left(3\,B-\tau\right)  \, P_{1,1,2}\\[4pt]
	&& \quad{}+2\, B^2 \,\left(3B-2\tau\right) \, P_{1,2,2}- 12\, A \,B^{2}\, P_{2,2,2} ,\\[8pt]
	h_2&=& 2AB_1\,\left(\,3BE  - 6CD - E\tau\,\right)  \\[4pt]
	&& \quad+ 12B_1^2\,C\,  P_{1,1,1}  - 2B_1^2 \,(3B  -\tau) \,P_{1,1,2} \\[4pt]
	&& \quad + 4 AB_1\,\left( 3B  - 2\tau\right)\, P_{1,2,2} - 24 \,A^2B_1 \,  P_{2,2,2}.      \end{array}\right.
\end{equation}
Note that $h_0=\det(\Hess')$, and that the cubic coefficients $P_{i,j,k}$ that appear here are the coefficients of the terms not involving $z$. For example, $P_{2,2,2}=\frac16H_{ppp}(0)$.

\begin{theorem}\label{thm:Type II NDC}
Consider an equilibrium point in $\RR^3$ with a Type II degeneracy; that is, $\tau\neq0$ and $h_0=0$. Then 
 \begin{enumerate}
 \item the vector field has a fold singularity provided $h_1,h_2$ do not both vanish; and
 \item in this case the family $H_\lambda=H-\lambda(\alpha q+\beta p+\gamma)$ (with $\alpha,\beta,\gamma\in\RR$) gives a versal unfolding of the singularity of the vector field, resulting in a saddle-node bifurcation of equilibria, provided
 $$(\beta,-\alpha,-\gamma)^T\; \not\in \;\Image(L).$$
 \end{enumerate}
 \end{theorem}

\begin{proof}
(i) Again, we rely on Lemma\,\ref{lemma:fold}. The key is that one needs to use different expressions for $\aa$ and $\vv$ depending on the values of $A,B,C$, and this leads to two separate non-degeneracy conditions: in fact it suffices to consider values of $B$ as follows. 

\smallskip

$$L=\begin{pmatrix}
B&2C&E\cr
-2A & -B+\tau  & -D\cr
0&0&\tau 
\end{pmatrix}.
$$

First, suppose $B\neq0$. Then we can use the non-zero vectors (recall $\tau\neq0$)
$$\aa=(2C,\;-B,\;0)^T,\quad \vv= -( 2A ,\; B ,\; \tfrac1\tau(B D - 2 A E)\,).$$
Then computing $\vv D^2(X_H)\aa^2$, after some simplification using $B(B-\tau)=4AC$, we find
$$\vv D^2(X_H)\aa^2=h_1$$
which for a fold we require to be non-zero in the case $B\neq0$ (note that $B$ is a factor of $h_1$). 

\smallskip
 Now suppose $B_1\neq0$ (that is, $B\neq\tau$). This time we use
$$\aa= (B_1,\, - 2 A , \, 0)^T   ,  \qquad \vv = \left(B_1,\, 2C, \, \tfrac1\tau(2CD - B_1E) \right)$$
With these vectors, both of which are non-zero, we obtain $\vv\, D^2(X_H)\,\aa^2=h_2$, which for a fold we require to be non-zero when $B_1\neq0$. (Note that $B_1$ is a factor of $h_2$.)

Since $\tau\neq0$,  $B$ and $B_1$ cannot both vanish simultaneously and part (i) of the theorem is proved,

\medskip
\noindent(ii) This is the statement of Theorem~\ref{thm:Type 2 fold}(ii) in this context. 
\end{proof}

%%%%%%%%%%%%%%%%%%%%%%%%%%%
\subsection{Type II saddle-node bifurcations}
We illustrate some cases of the theorem above in 3 dimensions. In the first example, we analyze the bifurcating equilibrium points, and in later examples we just record the condition on the 3-jet for the Type II equilibrium to be a fold. 

\begin{example} \label{eg:Type II}
Let $H = z - pq + p^2q - \lambda q$. At $\lambda=0$ this is a degenerate equilibrium  at the origin with eigenvalues $-1,-1,0$ and principal coefficient $\tau=-1$. We have $h_0=h_2=0$, but $h_1\neq0$ so the equilibrium is a fold singularity.   As $\lambda$ varies this family has a saddle-node bifurcation, with equilibria at $(q,p,z)=(0,\pm\sqrt\lambda,0)$ for $\lambda\geq0$. 

On one branch, the 0 eigenvalue becomes negative and on the other it becomes positive, as illustrated in  Figure\,\ref{fig:Type 2 eigenvalues}. 

\begin{itemize}
\item $(q,p,z) = (0,\sqrt{\lambda},0)$; at this point the eigenvalues are $-1,-1+2\sqrt\lambda$ and $-2\sqrt\lambda$ and the equilibrium is asymptotically stable (for small values of $\lambda$).
\item $(q,p,z) = (0,-\sqrt{\lambda},0)$; this point has one positive and two negative eigenvalues and so is unstable. 
\end{itemize}
\end{example}

There follows a table showing the non-degeneracy condition (up to a non-zero factor) for several simple values of $A,B,C$, and an admissible unfolding term. The unfolding term is independent of the values of $D,E$ and the $P_{i,j,k}$. In each of them, the analysis is similar to Example\,\ref{eg:Type II} above, and in fact that example is an instance of the penultimate of this list.
$$\begin{array}{cccc}
H_0 & h_1 & h_2  & \text{unfolding term} \\
\hline
A=B=C=0 & 0 &  P_{1,1,2} & p \\
A=B=0,C=1 & 0 & 6 P_{1,1,1} +\tau  P_{1,1,2} & p \\
A=1,B=C=0 & 0 & \tau  E +\tau^{2} P_{1,1,2} + 4 \tau  P_{1,2,2}+12 P_{2,2,2} & p \\
A=B_1=C=0 & E + P_{1, 2, 2} & 0 & q \\
A=1, B_1=C=0 &  \tau  (E +  P_{1,2,2})-6 P_{2,2,2} & 0 &p 
\end{array}
$$

%%%%%%%%%%%%%%%%%%%%%%%%%%%
%%%%%%%%%%%%%%%%%%%%%%%%%%%
\section{Legendre vector fields}
\def\Leg{\mathcal{L}}

In this section we consider the bifurcation theory of contact vector fields which are tangent to a given Legendre submanifold, and show the fact that they are the restriction of contact vector fields adds no restriction; that is, the bifurcation theory on the Legendre submanifold is the same as for generic vector fields in $\RR^n$. The main theorem is in essence due to Maschke \cite{Maschke-2016}.  Recall that a Legendre submanifold of $M^{2n+1}$ is a submanifold of dimension $n$ that is everywhere tangent to the contact structure, and that this is the maximal possible dimension of such a submanifold. 

In particular, we show in the theorem below that any given vector field (or family of vector fields) on a given Legendre submanifold can be extended to a contact vector field (or family of such) on the ambient contact manifold.  
 
Throughout this section we let $\Leg$ be a fixed Legendre submanifold of $(M,\xi)$, or of $\RR^{2n+1}$ as our analysis is local.  The following property of contact flows is due to Mruga{\l}a et al  \cite[Theorem 3]{Mrugala1991}, and is a direct analogue of a property of invariant Lagrangian submanifolds for symplectic flows. 

\begin{proposition}\label{prop:Legendre} Given a Hamiltonian $H$, the Legendre submanifold 
$\Leg$ is invariant under the flow of $X_H$ if and only if $H\restr{\Leg}\equiv0$.
\end{proposition}

\begin{proof}
Firstly if $\Leg$ is invariant under the flow, then at every point of $\Leg$ the vector field is tangent to $\Leg$ and hence is contained in the contact hyperplane, which implies $H=0$.

Conversely, suppose $H\restr{\Leg}=0$, and consider the flow induced by the vector field.  Now this flow preserves $H^{-1}(0)$ (as noted in \S\ref{sec:background}) and the vector field is therefore contained in the contact hyperplane at each point of $H^{-1}(0)$. If $\Leg$ is not invariant, let $x_0$ be a point where $X_H$ is not tangent to $\Leg$ and let $U$ be a neighbourhood of $x_0$ in $\Leg$ where this continues to hold. Consider the image of $U$ under the flow. This will be a submanifold of dimension $n+1$ tangent to the contact structure, which is not possible. 
\end{proof}

\paragraph{Generating functions} Following Arnold \cite[Appendix 4]{Arnold}, using Darboux coordinates one can  (locally) generate any Legendre submanifold of $\RR^{2n+1}$ as follows. 
Given a Legendre submanifold $\Leg\subset \RR^{2n+1}$, there is a subset  $I\subset\{1,\dots,n\}$ and a smooth function $S(q_i,p_a)$ ($i\in I, a\not\in I$) such that $\Leg$ is (locally) parametrized by $q_i,\,p_a$ by 
the following formulae
$$q_a=-\dbyd{S}{p_a},\quad p_i= \dbyd{S}{q_i},\quad z=S-p_a\frac{\partial S}{\partial p_a},$$
Conversely, given any such subest $I$ and generating function $S(q_i,p_a)$, the graph as given generates a Legendre submanifold.

As mentioned, the following is essentially due to Maschke \cite{Maschke-2016}, although there only for $S=S(q_i)$ (i.e., $I=\{1,\dots,n\}$), and without the (trivial) inclusion of parameters $\lambda$.

\begin{theorem}\label{thm:Legendre extension}
Let $\Leg$ be a Legendre submanifold of\/ $\RR^{2n+1}$ parametrized by $(q_i,p_a)$ as above, and let 
$$Y_\lambda =  f_j(q_i,p_a,\lambda)\dbyd{}{q_j} + f_b(q_i,p_a,\lambda)\dbyd{}{p_b}$$
be an arbitrary family of vector fields on $\Leg$, where the $f_i$ are smooth functions depending on parameter(s) $\lambda\in\RR^\ell$. Then there exists a family of contact vector fields $X_\lambda$ on a neighbourhood of $\Leg$ in $M$ whose restriction to $\Leg$ is equal to $Y_\lambda$: that is for $x\in\Leg$, $X_\lambda(x)=Y_\lambda(x)$. 
\end{theorem}

\begin{proof}
Let $\bar f_r(q,p,z,\lambda)=f_r(q_i,p_a,\lambda)$ and $\bar S(q,p,z)= S(q_i,p_a)$  be the trivial extensions of $f_r$ and $S$ respectively, to a  neighbourhood of $\Leg$ in $\RR^{2n+1}$ ($r=1,\dots,n$) --- that is, independent of $q_a,\,p_i,\,z$. 
Define
$$H(q,p,z,\lambda) = \left(p_i-\dbyd{\bar S}{q_i}\right)\bar f_i(q,p,z,\lambda) - \left(q_a+\dbyd{\bar S}{p_a}\right)\bar f_a(q,p,z,\lambda).$$
Clearly, $H\restr{\Leg}=0$ and hence $\Leg$ is invariant under the flow of $X_H$, which is to say, $X_H$ is tangent to $\Leg$. 
This implies that to check whether at points of $\Leg$ we have $X_H=Y$, we only need check the effect of $X_H$ on the coordinates $q_i,p_a$ of $\Leg$. Now, at points $x=(q,p,z)\in \Leg$, 
$$\begin{array}{lllllll}
X_H(q_i) &=&  \dot q_i &=& H_{p_i} &=& f_i(q_i,p_a,\lambda),\\
X_H(p_a) &=& \dot p_a &=& -H_{q_a}  -p_aH_z&=& f_a(q_i,p_a,\lambda),
\end{array}
$$
the latter since $H$ is independent of $z$. Hence the contact vector field $X_H$ coincides with $Y$ at points of $\Leg$, as required. 
\end{proof}

We remark that the extension chosen is in fact a conservative contact vector field, since $\Reeb(H)=H_z=0$. Had we allowed more general extensions of $f_i(q_i,p_a,\lambda)$ to $\bar f_i(q, p, z, \lambda)$ we would obtain other contact extensions of the vector field $Y$.

%%%%%%%%%%%%%%%%%%%%%%%%%%%
%%%%%%%%%%%%%%%%%%%%%%%%%%%
\setcounter{section}{0}
\appendix

\section{Recognizing fold singularities} \label{app:calculations}
\def\thesection{\Alph{section}}
In this appendix we derive a simple condition for recognizing when a map germ of corank 1 has a fold singularity, and when a deformation of a fold singularity is versal.  For details on $\kK$-equivalence see for example \cite{SBC} (note that $\kK$-equivalence is also called contact equivalence, but that could be confusing in the current context). 

Recall that for a map-germ $(\RR^n,0)\to(\RR^n,0)$, a \defn{fold} singularity is the least degenerate singularity, and is any germ $\kK$-equivalent to 
$$(x,\yy)\longmapsto(x^2,\yy)$$
with $x\in\RR$, $\yy\in\RR^{n-1}$. It has $\kK$-codimension 1, and 
$$(x,\yy;\lambda)\longmapsto(x^2-\lambda,\,\yy)$$
is a versal deformation (or unfolding).

\begin{lemma} \label{lemma:fold}
\begin{enumerate}
\item 
A corank-1 map-germ $F:(\RR^n,0)\to(\RR^n,0)$ has a fold singularity at the origin if and only if there are non-zero vectors $\aa\in\RR^3$ and $\vv\in\left(\RR^3\right)^*$ such that 
\begin{equation}\label{eq:conditions for fold}
\D F\,\aa=0,\quad \vv\, \D F = 0,\quad \vv\, \D^2F\,\aa^2\neq0,
\end{equation}
where the differentials are evaluated at the origin.
\item Given such a fold singularity, any 1-parameter deformation $F_\lambda$ (with $F_0=F$) is versal if and only if\, $\dot F(0) \not\in \Image(\d F(0))$, where $\dot F = \frac{\partial F_\lambda}{\partial \lambda}\restr{\lambda=0}$, and $\dot F(0)$ is its value at the origin.
\end{enumerate}
\end{lemma}

\begin{proof} Any map-germ $(\RR^n,0)\to(\RR^n,0)$ of corank 1 is contact equivalent to the map
\begin{equation}\label{eq:G corank 1}
G(x,\,\yy) =(g(x),\,\yy)
\end{equation}
for some smooth function-germ $g:(\RR,0)\to(\RR,0)$ with $g(0)=g'(0)=0$. Here $\yy\in\RR^{n-1}$ (see \cite[p.167]{SBC} for a proof.) 

(i)
Suppose the corank 1 map-germ  $F:(\RR^n,0)\to(\RR^n,0)$ is a fold singularity.  Then so is $G$ in \eqref{eq:G corank 1} and $g$ can be chosen to be $g(x)=x^2$ (and more generally $g''(0)\neq0$).  
Clearly, for $G$ the conditions \eqref{eq:conditions for fold} hold, with $\aa=(1,\; \zero)^T$ and $\vv=(1,\;\zero)$. Conversly, if $g''(0)=0$ then $G$ is not a fold and the condition $\vv D^2G\aa^2\neq0$. 

There remains to show that conditions \eqref{eq:conditions for fold} are unchanged under a contact equivalence.  This is a simple calculation, as follows.

Now $G$ and $F$ are $\kK$-equivalent iff $G(\xx)=A(\xx)F\circ\phi(\xx)$ where $A$ is an invertible $\xx$-dependent matrix and $\phi$ is a diffeomorphism  (all germs at 0). Now $G$ satisfies the conditions of the lemma: $\ww\, \D G=0,\;\D G\,\bb=0,\;\ww\, \D^2G\,\bb^2\neq0$.  Then 
$$\D G(\xx) = (\D A(\xx))F\circ\phi +  A(\xx)\D F(\phi(\xx))\D\phi(\xx)$$
so at $\xx=0$ where $F=G=0$ we have $\D G(0)=A(0)\,\D F(0)\,\D\phi(0)$. Then if $\ww \,\D G=0$ let $\vv=\ww\, A$, and if $\D G\,\aa=0$ let $\aa= \D\phi\,\bb$
then $\vv \,\D F=0$ and $\D F\,\aa=0$. Moreover,
\begin{equation*}
\begin{split}
\D ^2G = &(\D ^2A(\xx))F\circ\phi +  \D A(\xx)\D F(\phi(\xx))\D \phi(\xx)\\
& +A(\xx)\D ^2F(\phi(\xx))\D \phi(\xx)^2+A(\xx)\D F(\phi(\xx))\D ^2\phi(\xx)
\end{split}
\end{equation*}
Then putting $\xx=0$ and using $F(0)=G(0)=0$, we see
$$\ww\,\D^2G\,\bb^2= \vv \,\D^2F\,\aa^2$$
so that the latter is also nonzero.

(ii) A similar method of proof works here too: $F$ is $\kK$-equivalent to $G:(x,\yy)\mapsto(x^2,\yy)$ and the family $G_\lambda = G+\lambda\uu$ is a versal unfolding if and only if $\uu\not\in\RR\left\{e_2,\dots,\,e_n\right\}$, from the standard versality theorem for contact equivalence, eg.\ \cite{SBC}. 
\end{proof}

\begin{remark}\label{rmk:versal}
The importance of the deformations being \emph{versal} is that any two versal unfoldings of $\kK$-equivalent map-germs are themselves equivalent. In particular, any versal unfolding of a fold singularity is equivalent to the map $(x;\lambda)\mapsto x^2-\lambda$, showing that for $\lambda>0$ there are two zeros, while for $\lambda<0$ there are none (or vice versa if the sign of the $\lambda$ term is changed). That is, a versal unfolding of a fold is a saddle-node bifurcation.  However, $\kK$-equivalence of the vector fields does not respect eigenvalues. 
\end{remark}

\paragraph{Acknowledgements}
I would like to thank Alessandro Bravetti and Luis Garc\'ia-Naranjo for commenting on an early draft and suggesting some further references.

%%%%%%%%%%%%%%%%%%%%%%%%%%%%%%%%%%%%%%%

\printbibliography

%%%%%%%%%%%%%%%%%%%%%%%
\vskip 1cm

\noindent{ j.montaldi@manchester.ac.uk}

\bigskip

\noindent\begin{minipage}{0.3\textwidth}
\obeylines\it
Dept Mathematics
University of Manchester
Manchester M13 9PL
UK
\end{minipage}

\end{document}